\documentclass[preprint,12pt,review]{elsarticle}
\usepackage{epsfig}
\usepackage{amssymb}
\usepackage{amsthm}
\usepackage{lineno}
\journal{}
\begin{document}
\begin{frontmatter}
\title{TOPICS IN SPECIAL FUNCTIONS}
\author{Juan Carlos L\'opez Carre\~no.\corref{cor1}\fnref{label1}}
\ead{jclopez@unipamplona.edu.co}
%% \ead[url]{home page}
\fntext[label1]{Universidad de Pamplona, Departamento de Matem\'aticas}
\author{Rosalba Mendoza Su\'arez.\corref{cor1}\fnref{label1}}
\ead{rosalbame@unipamplona.edu.co}
\author{Jairo Alonso Mendoza S.\corref{cor1}\fnref{label2}}
\ead{jairoam@unipamplona.edu.co}
\fntext[label2]{Universidad de Pamplona, Departamento de F\'isica y Geolog\'ia}
\address{$^{1}$Departamento de Matem\'aticas, Universidad de Pamplona \\
Pamplona, Norte de Santander, Colombia.}
\address{$^{2}$Departamento de F\'isica y Geolog\'ia, Universidad de Pamplona \\
Pamplona, Norte de Santander, Colombia. Grupo de investigacion INTEGRAR}
\begin{abstract}
In this article, we solve the connection problem of the Hermite polynomials
with the classical continuous orthogonal polynomials belonging to Askey
scheme, using the hypergeometric functions method combined is with the
work the Fields and Ismail.
\end{abstract}
\end{frontmatter}
\section{INTRODUCTION}
\label{sec:introduction}
\noindent The connection problem is to find the coefficients
$c_{nk}$ in the expansion of one polynomial $P_n(x)$ in terms of an
arbitrary sequence of orthogonal polynomials $\{Q_k(x)\}$,
\begin{equation}\label{1}
P_n(x)=\sum_{k=0}^{n}c_{nk}Q_k(x).
\end{equation}
\noindent A wide variety of methods have been devised for computing
the connection coefficients $c_{nk}$, either in closed form or by
means of recurrence relations, usually in $k$. Lewanowicz
\cite{Lewa} has shown that the connection problem (\ref{1}) can,
sometimes be solved by taking advantage of known results from the
theory of generalized hypergeometric functions, derived by Fields
and Wimp \cite{fieldswimp}.
\vskip 0.5cm

\noindent The hypergeometric functions method has been devised to
solve the connection problem, that involving classical orthogonal
polynomials . \cite{Lewa,R.Askey,luk2,ismail,fielismail,Jorge5,Jorge4}.
These methods have also been used by S\'anchez-Ruiz \cite{Jorge5} to
obtain the connection formulae involving squares of Gegenbauer
Polynomials. In \cite{Jorge4}, authors seeking for solving (\ref{1})
for a much wider class of polynomials,  defined by terminating
hypergeometric series, obtain connection formulae for Wilson and
Racah polynomials with special parameter values. They also solve the
connection problem for the families of generalized Jacobi and
Laguerre polynomials defined by Sister Celine.
\vskip 0.5cm

\noindent In \cite {ismail} authors consider the expansion of
arbitrary power series in various classes of polynomial sets. In
particular,  they obtain the connection formulae which generalize
the expansions formulae Verma and Field and Wimp. (See
Lema 2.1, Lema 2.2).

\section{NOTATION AND PRELIMINARY RESULTS}
The generalized hypergeometric function is defined by
\begin{equation}\label{hyp-def}
_pF_q \Bigg( \begin{array}{c} {a_1,a_2,\ldots,a_p} \\
{b_1,b_2,\ldots,b_q}
\end{array}
\Bigg|x\Bigg)=\displaystyle \sum_{k=0}^\infty
\frac{(a_1)_k(a_2)_k\cdots (a_p)_k}{(b_1)_k(b_2)_k\cdots
(b_q)_k}\frac{x^k}{k!},
\end{equation}
where $(a)_n$ represents the Pochhammer symbol and it used in the
theory of special functions to represent the rising
factorial.
\begin{eqnarray*}
(a)_n := a(a+1)(a+2)\cdots(a+n-1)= \frac{\Gamma (a+n)}{\Gamma (a)}=\frac{%
(-1)^n\,\Gamma (1-a)}{\Gamma (1-a-n)},
\end{eqnarray*}
\begin{eqnarray*}
(a)_0= 1,
\end{eqnarray*}
where $a_i \in C,$ $1\leq i \leq p,$ $b_j\in C,$ $1\leq j\leq q,$
with $b_j  \in \hskip-0.3cm/ \hskip0.01cm N_0.$
\noindent Throughout this article, the letters $p,q,r,s,t,u$ and $n$
stand for nonnegative integers.
\noindent We call $x$ the argument of the function, and $a_j$, $b_j$
the parameters.
\noindent To shorten the notation for the left-hand side of
(\ref{hyp-def}),   we will write it as
\begin{equation}
_pF_q \Bigg(
\begin{array}{c} {[a_p]}\\ {[b_q]}
\end{array} \Bigg| x \Bigg) = \sum_{k=0}^\infty \, \frac{[a_p]_{k} \, x^k}
{[b_q]_k \, k!} \,,
\end{equation}
where $[a_p]$ and $[b_q]$ represent the sets $\{a_1, a_2,...,a_p \}$
and $\{b_1, b_2,...,b_q \}$, respectively. We use the abbreviated
notation
\begin{equation}
{[a_p]_k} = \prod_{i=1}^p (a_i)_k, \, \quad  \quad {[b_q]_k} =
\prod_{j=1}^q (b_j)_k.
\end{equation}
To prove the theorems in section 3, we use known results
from the theory of generalized hypergeometric
functions.\\
\noindent Lemma 2.1 (See \cite{ismail}, Formulae 1.3
and 3.2)
\begin{equation} \label{1.3}
\begin{array}{rl}
\sum_{m=0}^{\infty} \,a_m b_m \frac{(zw)^m}{m!}  = &\sum_{n=0}^{\infty}
\,\frac{(-z)^n}{n!(\gamma+n)_n}  \sum_{r=0}^{\infty}
\frac{(\mu)_{n+r} (\theta)_{n+r}}{r!(\gamma+2n+1)_r} b_{n+r} \ z^r \\
 &  \\
& \times  \sum_{s=0}^{n}
\,\frac{(-n)_s}{s!}\frac{(n+\gamma)_s}{(\mu)_s(\theta)_s}a_sw^s.
\end{array}
\end{equation}
\begin{eqnarray} \label{3.2}
\sum_{m=0}^{\infty} \,a_m b_m (zw)^m= \sum_{n=0}^{\infty}
\,\frac{(c)_n(-z)^n}{n!} \sum_{j=0}^{\infty} \, \frac{(n+c)_{j}
}{j!} b_{n+j} \ z^j\, \sum_{k=0}^{n} \,\frac{(-n)_k}{(c)_k}a_kw^k.
\end{eqnarray}
\noindent Lemma 2.2  (See \cite{fieldswimp}, \cite[Vol. II, p.
7]{luk2})
\begin{equation}\label{luke}
\begin{array}{rl}
_{p+r}F_{q+s} \bigg(\begin{array}{c} {[a_p]},{[c_r]} \\
{[b_q]},{[d_s]}
\end{array} \bigg| z\omega\bigg)= & \sum_{n=0}^\infty
\frac{(a_p)_n(c)_n(-z)^n}{(b_q)_nn!}
\, _{p+1}F_{q} \bigg(\begin{array}{c} {[n+\alpha]},{[n+a_p]} \\
{[n+b_q]}
\end{array} \bigg| z\bigg) \\
&\times  _{r+1}F_{s+1} \bigg(\begin{array}{c} -n,{[c_r]} \\
c,{[d_s]}
\end{array} \bigg| \omega\bigg).
\end{array}
\end{equation}
\begin{eqnarray}\label{wimp2}
\begin{array}{rl}
&_{p+r+1}F_{q+s} \Bigg(
\begin{array}{c} -n, {[a_p]}, {[c_r]} \\ {[b_q], {[d_s]}}
\end{array} \Bigg| zw \Bigg) =  \sum_{k=0}^{n} \, \Bigg( \begin{array}{c}  n \\ k  \end{array}\Bigg) \,
\frac{[a_p]_k \, [\alpha_t]_k \, z^k}{[b_q]_k[\beta_u]_k} \\
& \times \, _{p+t+1}F_{q+u} \Bigg(
\begin{array}{c} k-n, {[k+a_p]}, {[k+\alpha_t]}\\ {[k+b_q]}, {[k+\beta_u]}
\end{array} \Bigg| z \Bigg) \, \,  _{r+u+1}F_{s+t} \Bigg(
\begin{array}{c} -k, {[c_r]}, {[\beta_u]}\\ {[d_s], {[\alpha_t]}}
\end{array} \Bigg| w \Bigg).
\end{array}
\end{eqnarray}
\noindent Lemma 2.3 \, \,  From the generalized hypergeometric
definition, it is obtained:
\begin{equation}\label{esteban}
\begin{array}{rl}
&_pF_q \Bigg(
\begin{array}{c} {[a_p]}\\ {[b_q]}
\end{array} \Bigg| z \Bigg) =
\, _{2p}F_{2q+1} \Bigg(\begin{array}{c} \left
[\frac{a_p}{2}\right],\left[\frac{a_p+1}{2}\right]\\\frac{1}{2},\left[\frac{b_q}{2}\right],\left[\frac{b_q+1}{2}\right]
\end{array} \bigg| 4^{p-q-1}z^2\Bigg) +\frac{(\prod_{i=1}^pa_i)z}{(\prod_{i=1}^q b_i)} \\
& _{2p}F_{2q+1} \Bigg(\begin{array}{c}
\left[\frac{a_p+1}{2}\right],\left[\frac{a_p+2}{2}\right]\\\frac{3}{2},\left[\frac{b_q+1}{2}\right],\left[\frac{b_q+2}{2}\right]
\end{array} \bigg| 4^{p-q-1}z^2\Bigg).\\
\end{array}
\end{equation}
\section{Results}
\noindent In this section, we use the lemmas 2.1,
2.2, 2.3 and the different hypergeometric series representation the classical continuous
orthogonal polynomials to solve the connection problem among them.

\subsection{Theorem  (\textrm{Connection formula for Laguerre Polynomials in series of Hermite Polynomials})}

\begin{eqnarray*}
L_n(x)= \, \sum_{k=0}^{n} \,\frac{(-n)_{k}} {2^{k}(k)!(k)!}\,
_{2}F_{2} \bigg(\begin{array}{c} \frac{k-n}{2},  \frac{k+1-n}{2}\\
\frac{k+1}{2}, \frac{k+2}{2},
\end{array} \bigg| \frac{1}{4} \bigg) \,\,H_{k}(x).
\end{eqnarray*}
\textbf{Proof:} Some of the hypergeometric representations of
Laguerre and Hermite polynomials are given by:
\begin{equation}\label{30}
L_n(x)= \sum_{k=0}^n \frac{(-n)_k}{k!}\frac{x^k}{k!}=\, _1F_1
\bigg(\begin{array}{c} -n \\ 1
\end{array} \bigg| x \bigg),
\end{equation}
\begin{equation}\label{31}
H_{2m}(x)=(-1)^m2^{2m}\left(\frac{1}{2}\right)_m \, _1F_1
\bigg(\begin{array}{c} -m \\ \frac{1}{2}
\end{array} \bigg| x^2 \bigg),
\end{equation}
\begin{equation}\label{32}
H_{2m+1}(x)=(-1)^m2^{2m+1}\left(\frac{3}{2}\right)_m \,x \,  _1F_1
\bigg(\begin{array}{c} -m \\ \frac{3}{2}
\end{array} \bigg| x^2 \bigg).
\end{equation}
\noindent Using formula (\ref{esteban}) of Lemma 2.3 with
the following identification
$$p=1, \, \, q=1, \, \, \{a_1\}=\{-n\}, \, \,\{b_1\}=\{1\}, \,\, z=x $$
it is obtained
\begin{equation}\label{33}
_{1}F_{1} \bigg(\begin{array}{c} -n \\1
\end{array} \bigg| x\bigg)= \, _{2}F_{3} \bigg(\begin{array}{c} -\frac{n}{2}, \frac{-n+1}{2}\\
\frac{1}{2}, \frac{1}{2}, 1
\end{array} \bigg| \frac{1}{4}x^2 \bigg) \,\, +x(-n)\,  _{2}F_{3} \bigg(\begin{array}{c} \frac{-n+1}{2},
\frac{-n+2}{2}\\\frac{3}{2}, 1, \frac{3}{2}
\end{array} \bigg| \frac{1}{4}x^2\bigg).
\end{equation}
If we call:
\begin{equation}\label{36}
A_1(x)=\,  _{2}F_{3} \bigg(\begin{array}{c} -m, \frac{-2m+1}{2}\\
\frac{1}{2}, \frac{1}{2}, 1
\end{array} \bigg| \frac{1}{4}x^2 \bigg),
\end{equation}
\begin{equation}\label{89}
A_2(x)= -2mx \,  _{2}F_{3} \bigg(\begin{array}{c} -(m-1),
\frac{-2m+1}{2}\\\frac{3}{2}, 1, \frac{3}{2}
\end{array} \bigg| \frac{1}{4}x^2\bigg).
\end{equation}
and $n=2m$,  then  (\ref{33}) can be written as:
\begin{equation}\label{35}
L_n(x) = A_1(x)+ A_2(x).
\end{equation}
Using Fields-Wimp formula  (\ref{luke}) we connect $A_1(x)$ give in
(\ref{36}) with the Hermite polynomials of even degree, (\ref{31}).
making the following identification:
$$ [d_s]=\emptyset, \, \, s=0, \, \, [c_r]=\emptyset, \, \, r=0, \,c=\frac{1}{2} \, [a_p]=\left\{-m,\frac{-2m+1}{2}\right\} $$
$$p=2, \, \, [b_q]= \{\frac{1}{2}, \frac{1}{2}, 1\}, \, \, q=3, \, \, w=x^2, \, \, z=\frac{1}{4}, $$
it is obtained
\begin{equation}\label{34}
\begin{array}{rl}
_{2}F_{3} \bigg(\begin{array}{c} -m, \frac{-2m+1}{2}
\\\frac{1}{2}, \frac{1}{2}, 1
\end{array} \bigg| \frac{x^2}{4}\bigg)=& \, \sum_{k=0}^m  \,\frac{(-m)_k\left(\frac{-2m+1}{2}\right)_k\left(\frac{1}{2}\right)_k\left(\frac{-1}{4}\right)^k}
{\left(\frac{1}{2}\right)_k\left(\frac{1}{2}\right)_k\left(1\right)_kk!} \\  & _{3}F_{3} \bigg(\begin{array}{c} k+\frac{1}{2}, k-m, k+\frac{1-2m}{2}\\
k+\frac{1}{2}, k+\frac{1}{2},k+ 1
\end{array} \bigg| \frac{1}{4} \bigg) \,\, _{1}F_{1} \bigg(\begin{array}{c}
-k\\\frac{1}{2}\end{array} \bigg| x^2\bigg).
\end{array}
\end{equation}
The left side of  (\ref{34}) is $A_1(x)$ given in (\ref{36}) by
reorganizing the terms of the right side of  (\ref{34})the  Hermite
polynomials of even degree arise(\ref{31}),and using the formula for
duplicating the  Pochhammer symbol, it results that:
\begin{eqnarray}\label{39}
A_1(x)= \, \sum_{k=0}^m \,\frac{(-2m)_{2k}}
{2^{2k}(2k)!(2k)!}\, _{2}F_{2} \bigg(\begin{array}{c} \frac{2k-2m}{2},  \frac{2k+1-2m}{2}\\
\frac{2k+1}{2}, \frac{2k+2}{2},
\end{array} \bigg| \frac{1}{4} \bigg) \,\,H_{2k}(x).
\end{eqnarray}
Likewise, polynomials $A_2(x)$ of  (\ref{89}) are connected with
Hermite polynomials of odd degree, (\ref{32}); making the following
identification:
$$[d_s]=\emptyset, \, s=0, \, [c_r]=\emptyset, \, r=0, c=\frac{3}{2}, [a_p]=\left\{-\frac{(2m-2)}{2}, \frac{-2m+1}{2}\right\},
\,p=2$$
$$ [b_q]=\left\{\frac{3}{2}, \frac{3}{2}, 1\right\}, \, q=3, w=x^2, \, z=\frac{1}{4},$$
\begin{equation}\label{38}
A_2(x)= \, \sum_{k=0}^{m-1} \,\frac{(-2m)_{2k+1}}
{2^{2k+1}(2k+1)!(2k+1)!}\,
_{2}F_{2} \bigg(\begin{array}{c} \frac{2k+1-2m}{2},  \frac{2k+1-2m+1}{2}\\
\frac{2k+1+2}{2}, \frac{2k+1+1}{2},
\end{array} \bigg| \frac{1}{4} \bigg)
\,\,H_{2k+1}(x).
\end{equation}
Recalling $n=2m,$ and $L_n(x)=A_1(x)+A_2(x)$  then from (\ref{39})
and  (\ref{38}) we conclude that
\begin{eqnarray*}
L_n(x)= \, \sum_{k=0}^{n} \,\frac{(-n)_{k}} {2^{k}(k)!(k)!}\,
_{2}F_{2} \bigg(\begin{array}{c} \frac{k-n}{2},  \frac{k+1-n}{2}\\
\frac{k+1}{2}, \frac{k+2}{2},
\end{array} \bigg| \frac{1}{4} \bigg) \,\,H_{k}(x).
\end{eqnarray*}
\subsection{Theorem (Connection formula for Hermite polynomials in series of Laguerre Polynomials)}
\begin{eqnarray*}
H_N(x)= N!2^N \sum_{m=0}^{N} \, _{2}F_{2} \bigg(\begin{array}{c} -\frac{1}{2}(N-m), -\frac{1}{2}(N-m-1)\\
-\frac{1}{2}(N), -\frac{1}{2}(N-1)\end{array} \bigg| -\frac{1}{4}
\bigg)\frac{(-N)_{m}}{(m)!}L_{m}(x),
\end{eqnarray*}
\textbf{Proof:}  Hermite polynomials can be expressed as:
\begin{eqnarray*}
H_N(x)= N! \sum_{k=0}^{[\frac{N}{2}]} \,
\frac{(-1)^k}{k!(N-2k)!}(2x)^{N-2k},
\end{eqnarray*}
\noindent If in the above expression we set up $N=2p$, this is
rewritten as:
\begin{equation}{\label{hermite2}}
H_{2p}(x)= (2p)! \sum_{k=0}^p \,
\frac{(-1)^k}{k!(2p-2k)!}(2x)^{2p-2k}.
\end{equation}
\noindent Let $m:=2p-2k$, for  $k=0,1,2,...,p$ one realized that $m$
goes through the set $\{2p,2p-2,...,0\}$
\noindent If in (\ref{hermite2}) we make the following substitution
$m=2p-2k$, it is obtained:
\begin{equation}\label{58}
H_{2p}(x)= (2p)! \sum_{m\in \{2p, 2p-2,...,0\}} \,
\left\{\frac{(-1)^{\frac{2p-m}{2}}2^mx^m}{m!(\frac{2p-m}{2})!}\right\}
\end{equation}
\noindent With the following identification   $a_k=\frac{1}{k!}$,
$c=1$, $w=x,$ $a_m=\frac{1}{m!}$ $z=1$ and
\begin{eqnarray}\label{bm}
b_m =\Bigg\{ \begin{array}{c}\frac{(-1)^{\frac{2p-m}{2}}2^m(2p)!}{(\frac{2p-m}{2})!},
m \in \{2p, 2p-2,..., 0\}  \\
0, \hskip0.5cm   \textrm{otherwise}
\end{array}
\end{eqnarray}
\noindent the left side of formula  (\ref{3.2}) of Lemma 2.1,is
precisely the expression for Hermite polynomials given in
(\ref{58}). By reorganizing the terms of the right side of
(\ref{3.2}) Laguerre polynomials result.
\begin{eqnarray*}
H_{2p}(x)= \sum_{n=0}^{\infty} \, (-1)^n\sum_{j=0}^{\infty} \,
\frac{(n+1)_j}{j!}b_{n+j}L_n(x)
\end{eqnarray*}
Keeping in mind the definition given for coefficients $b_m$,
 $H_{2p}(x)$ it can be written:
\begin{eqnarray*}
\begin{array}{rl}
H_{2p}(x)= & \left\{\sum_{k=0}^{p} \,\sum_{s=0}^{p-k} \,(-1)^{2k}
\frac{(2k+1)_{2s}}{(2s)!}b_{2k+2s}L_{2k}(x)\right\} \\
& + \left\{\sum_{k=0}^{p-1} \, \sum_{s=0}^{p-k-1} \,(-1)^{2k+1}
\frac{(2k+1+1)_{2s+1}}{(2s+1)!}b_{2k+1+2s+1}L_{2k+1}(x)\right\}.
\end{array}
\end{eqnarray*}
\noindent Let
\begin{eqnarray*}
H_{2p}(x):= I + II.
\end{eqnarray*}
Making use of the expression of $b_m$ given in (\ref{bm}), we have
\begin{eqnarray*}
I=\sum_{k=0}^{p} \,\sum_{s=0}^{p-k} \,(-1)^{2k}
\frac{(2k+1)_{2s}}{(2s)!}
\frac{(-1)^{p-k-s}2^{2k+2s}(2p)!}{(p-k-s)!}\, L_{2k}(x).
\end{eqnarray*}
\begin{eqnarray*}
II=\sum_{k=0}^{p-1} \, \sum_{s=0}^{p-k-1} \,(-1)^{2k+1}
\frac{(2k+2)_{2s+1}}{(2s+1)!}
\frac{(-1)^{p-k-s-1}2^{2k+2s+2}(2p)!}{(p-k-s-1)!}\, L_{2k+1}(x).
\end{eqnarray*}
By distributing some terms, the expression of  $I$ can be written
as:
\begin{equation}\label{I}
I= (2p)!2^{2p} \sum_{k=0}^{p} \,\left\{\sum_{s=0}^{p-k} \,
\frac{(2k+1)_{2s}(-1)^{p-k-s}2^{2k+2s}}
{(2s)!(p-k-s)!}\frac{(2k)!}{2^{2p}(-2p)_{2k}}\right\} \, \,
\frac{(-2p)_{2k}}{(2k)!}L_{2k}(x).
\end{equation}
\noindent Using some properties of the Pochhammer symbol, the
interior sum of (\ref{I}), performed over $s$, can be written as:
\begin{eqnarray}\label{11}
& & \sum_{s=0}^{p-k} \,
\frac{(k-p+\frac{1}{2})_{p-k-s}}{(\frac{1}{2}-p)_{p-k-s}}\frac{(k+s)!(p-k)!
(-1)^{p-k-s}2^{2k+2s}}{s!
p!(p-k-s)!2^{2p}}\frac{2^{2s}(-\frac{1}{4})^s}{(-1)^s}.
\end{eqnarray}
\noindent With the substitution
$$2k+2s=2p-2t,$$
\noindent the above expression can be written as:
\begin{eqnarray*}
\sum_{t=0}^{p-k} \,
\frac{(-(p-k))_t(k-p+\frac{1}{2})_{t}(-1)^t}{(-p)_tt!(\frac{1}{2}-p)_{t}2^{2t}}.
\end{eqnarray*}
As $N=2p$, then(\ref{I})is transformed into:
\begin{equation}\label{71}
I= N!2^N \sum_{k=0}^{\frac{N}{2}} \, _{2}F_{2} \bigg(\begin{array}{c} -\frac{1}{2}(N-2k), -\frac{1}{2}(N-2k-1)\\
-\frac{1}{2}(N), -\frac{1}{2}(N-1)
\end{array} \bigg| -\frac{1}{4} \bigg)
\, \frac{(-N)_{2k}}{(2k)!}L_{2k}(x).
\end{equation}
Likewise  $II$, turns out to be
\begin{equation}\label{72}
II= N!2^N \sum_{k=0}^{\frac{N}{2}-1} \, _{2}F_{2} \bigg(\begin{array}{c} -\frac{1}{2}(N-2k-1), -\frac{1}{2}(N-2k-2)\\
-\frac{1}{2}(N), -\frac{1}{2}(N-1)
\end{array} \bigg| -\frac{1}{4} \bigg)
\, \frac{(-N)_{2k+1}}{(2k+1)!}L_{2k+1}(x).
\end{equation}
\noindent Recalling that
\begin{eqnarray*}
H_N(x)= I+II.
\end{eqnarray*}
Thence,  (\ref{71}) and (\ref{72}) allow us to write the connection
of Hermite $H_N(x)$ with Laguerre polynomials:
\begin{eqnarray*}
& & H_N(x)= N!2^N \sum_{k=0}^{\frac{N}{2}} \, _{2}F_{2} \bigg(\begin{array}{c} -\frac{1}{2}(N-2k), -\frac{1}{2}(N-2k-1)\\
-\frac{1}{2}(N), -\frac{1}{2}(N-1)
\end{array} \bigg| -\frac{1}{4} \bigg)\, \frac{(-N)_{2k}}{(2k)!}L_{2k}(x)\nonumber\\
& & \nonumber\\
& & + N!2^N \sum_{k=0}^{\frac{N}{2}-1} \, _{2}F_{2} \bigg(\begin{array}{c} -\frac{1}{2}(N-2k-1), -\frac{1}{2}(N-2k-2)\\
-\frac{1}{2}(N), -\frac{1}{2}(N-1)
\end{array} \bigg| -\frac{1}{4} \bigg) \frac{(-N)_{2k+1}}{(2k+1)!}L_{2k+1}(x)
\end{eqnarray*}
\begin{eqnarray*}
H_N(x)= N!2^N \sum_{m=0}^{N} \, _{2}F_{2} \bigg(\begin{array}{c} -\frac{1}{2}(N-m), -\frac{1}{2}(N-m-1)\\
-\frac{1}{2}(N), -\frac{1}{2}(N-1)\end{array} \bigg| -\frac{1}{4}
\bigg)\frac{(-N)_{m}}{(m+)!}L_{m}(x),
\end{eqnarray*}
\noindent \textbf{\emph{Note:}} \, Likewise the procedure can
be applied when $N$ is odd.
\subsection{ Theorem (Connection formula for Hermite polynomials in series of Shifted Jacobi Polynomials)}
\begin{eqnarray*}
\begin{array}{rl}
H_n(x)= & \frac{(-n)_m4^n(2m+\lambda)(\alpha+1)_n}{(\alpha+1)_m(\lambda+m)_{n+1}} \, _{4}F_{2} \bigg(\begin{array}{c} \Delta(2,m-n), \Delta(2,-\lambda-n-m)\\
\Delta(2,-\alpha-n)\end{array} \bigg| \frac{1}{4}
\bigg) \\
&\times \ P_{m}^{(\alpha,\beta)}(1-x)
\end{array}
\end{eqnarray*}
where, $\lambda=\alpha+\beta+1,$ y
$\Delta(r;\varphi)=\frac{(\varphi+j-1)}{r} , \, \,
j=1,\cdot\cdot\cdot,r.$\\
\noindent \textbf{Proof:}  Shifted Jacobi polynomials can be written
as:
\begin{eqnarray*}
R_n^{(\alpha,\beta)}(x)= \frac{(-1)^n(\beta+1)_n}{n!}\, _{2}F_{1}
\bigg(\begin{array}{c} -n, n+\alpha+\beta+1\\\beta+1\end{array}
\bigg|x\bigg)
\end{eqnarray*}
\begin{eqnarray*}
R_n^{(\alpha,\beta)}(x)= \frac{(-1)^n(\beta+1)_n}{n!}\,\sum_{s=0}^n
\, \frac{(-n)_s(n+\alpha+\beta+1)_s}{(\beta+1)_s}\frac{x^s}{s!}
\end{eqnarray*}
\noindent With the following substitution $\gamma=\alpha+\beta+1,$
$\theta=\beta+1,$ $\mu=1,$ $a_s=s!$,$w=x,$ $z=1$ and
\begin{eqnarray}\label{bm}
b_m =\Bigg\{
\begin{array}{c}\frac{(-1)^{\frac{2p-m}{2}}2^m(2p)!}{m!(\frac{2p-m}{2})!},
\, \, m \in \{2p, 2p-2,..., 0\}\\
0 \hskip0.5cm  \textrm{otherwise}
\end{array}
\end{eqnarray}
\noindent the left side of formula (\ref{1.3}) of Lemma 2.1, is
precisely the expression for Hermite polynomials given in
(\ref{58}). On the other hand in the right side  (\ref{1.3}) the
Shifted Jacobi, polynomials,appear. Therefore (\ref{1.3}) is
transformed in
\begin{eqnarray*}
H_{2p}(x)= \sum_{n=0}^{\infty} \,
\frac{1}{(\beta+1)_n(\alpha+\beta+1)_n}\sum_{r=0}^{\infty} \,
\frac{(n+r)!(\beta+1)_{n+r}}{r!(\alpha+\beta+2n+2)_r}b_{n+r}R_n^{(\alpha,\beta)}(x).
\end{eqnarray*}
The result is also obtained by a similar procedure to the one
followed in the proof of Theorem 3.2.
\subsection{Theorem (Connection formula for Shifted Jacobi
Polynomials  in series of Hermite polynomials)}
\begin{eqnarray*}
& & R_n^{(\alpha,\beta)}(x)=\sum_{j=0}^{n} \,
\frac{(-1)^n(-1)^j(\beta+1)_n (n+\alpha+\beta+1)_j }{(n-j)! 2^{j} \,
(j)!
\,(\beta+1)_{j}}\nonumber \\
& & \times \, _{4}F_{2} \bigg(
\begin{array}{c} \frac{j-n}{2}, \frac{j+1-n}{2}, \frac{j+n+\alpha+\beta+1}{2},
\frac{j+n+\alpha+\beta+2}{2} \\\frac{j+\beta+1}{2},
\frac{j+\beta+2}{2}
\end{array} \bigg| 1\bigg) \,H_{j}(x).
\end{eqnarray*}
\textbf{Proof:} Using formula(\ref{esteban}) of Lemma 2.3 with the
following identification
$$p=2 \quad q=1 \quad a_1 =\{-n\} \quad a_2 =\{n+\alpha+\beta+1\} \quad b_1 =\{\beta+1\}  z= x $$
it is obtained
\begin{eqnarray}\label{rosalba}
& &R_n^{(\alpha,\beta)}(x) = \, (-1)^n\frac{(\beta+1)_n}{n!} \,_4F_3 \bigg(\begin{array}{c} -\frac{n}{2}, \frac{n+\alpha+\beta+1}{2},\frac{1-n}{2},\frac{n+\alpha+\beta+2}{2}\\
\frac{1}{2}, \frac{\beta+1}{2}, \frac{\beta+2}{2}\end{array} \bigg| x^{2} \bigg)\nonumber \\
& & +(-1)^n\frac{(\beta+1)_n}{n!} \,\frac{x
(-n)(n+\alpha+\beta+1)}{\beta+1} \, \,
_4F_3 \bigg(\begin{array}{c} \frac{1-n}{2}, \frac{n+\alpha+\beta+2}{2},\frac{2-n}{2},\frac{n+\alpha+\beta+3}{2}  \\
\frac{3}{2},\frac{\beta+2}{2}, \frac{\beta+3}{2}\end{array} \bigg|
x^2 \bigg) \nonumber \\ 
\end{eqnarray}
If we note by
\begin{eqnarray*}
A_1(x)= \, (-1)^n\frac{(\beta+1)_n}{n!} \,_4F_3 \bigg(\begin{array}{c} -\frac{n}{2}, \frac{n+\alpha+\beta+1}{2},\frac{1-n}{2},\frac{n+\alpha+\beta+2}{2}\\
\frac{1}{2}, \frac{\beta+1}{2}, \frac{\beta+2}{2}\end{array} \bigg|
x^{2} \bigg).
\end{eqnarray*}
\begin{eqnarray*}
A_2(x)= (-1)^n\frac{(\beta+1)_n}{n!} \,\frac{x
(-n)(n+\alpha+\beta+1)}{\beta+1} \, \,
_4F_3 \bigg(\begin{array}{c} \frac{1-n}{2}, \frac{n+\alpha+\beta+2}{2},\frac{2-n}{2},\frac{n+\alpha+\beta+3}{2}  \\
\frac{3}{2},\frac{\beta+2}{2}, \frac{\beta+3}{2}\end{array} \bigg|
x^2 \bigg)
\end{eqnarray*}
then  (\ref{rosalba})  can be written as
\begin{eqnarray*}
R_n^{(\alpha,\beta)}(x)= A_1 (x) + A_2(x).
\end{eqnarray*}
Following a similar procedure to that used in  Theorem 3.1 one end
up with the expected result.

\section{Summary and continuity of the work}
In this work, we present the connections between the Polynomials of Laguerre in series of the polynomials of Hermite, between the polynomials of Hermite and the polynomials of Laguerre, between the polynomials of Hermite and the polynomials of Shifted Jacobi, between the polynomials of Shifted Jacobi and the polynomials of Hermite. He expects to present it in a new work the connection between the polynomials of Bessel and the polynomials of Hermite, between the polynomials of Hermite and the polynomials of Bessel.

\noindent {\bf{References}}

\bibliographystyle{elsarticle-harv}

\begin{thebibliography}{00}

\bibitem[S.Lewanowicz (1997)]{Lewa} S. Lewanowicz, {\em The hypergeometric functions approach to the connection
problem for the classical orthogonal polynomials}, Preprint (1997).


\bibitem[J.L.Fields(1961)]{fieldswimp} J. L. Fields, J. Wimp, {\em Expansions of hypergeometric functions
in hypergeometric functions}, Math. Comp. {\bf 15} (1961), 390--395.


\bibitem [R. Askey(1975)]{R.Askey} R. Askey, {\em Orthogonal Polynomials and Special Functions},
Regional Conference Series in Applied Mathematics, Volume 21 (SIAM,
Philadelphia, 1975).


\bibitem[Y.L. Luke (1969)]{luk2} Y. L. Luke, {\em The Special Functions and their Approximations}
(Academic, New York, 1969).

\bibitem[Jerry L. (1975)]{ismail}Jerry L. Fields and Mourad E. H. Ismail. {\em Polynomial expansion,} Math. Comp. {\bf 131} (1975), 894-902.

\bibitem[Jerry L. (196)]{fielismail} Jerry L. Fields and J:Wimp, {/em  Expansions of hypergeometric functions in    hypergeometric functions,} Math. Comp. {\bf 15} (1961), 390-395.

\bibitem[J. Sanchez (2001)]{Jorge5} J. S\'anchez-Ruiz, {\em Linearization and connection formulae involving squares of Gegenbauer  polynomials}, Appl. Math. Lett. {\bf 14} (2001), 261--267.

\bibitem[J. Sanchez 2 (2001)]{Jorge4} J. S\'anchez-Ruiz, J. S. Dehesa, {\em Some connection and linearization problems for polynomials in and beyond the Askey scheme}, J. Comput. Appl. Math. {\bf 133} (2001), 579--591.

\bibitem[J. Sanchez 3 (2003)]{Jorge3} S\'anchez R J, P. L. Artés, A. Mart\'inez-Finkelshtein, J. S. Dehesa, {\em Expansions in series of varing Laguerre polynomials and some applications to molecular potentials}, J.Comput.  Appl.Math.153 (2003), 411-421.
\end{thebibliography}

\end{document}